\documentclass[12pt]{amsart}
\usepackage{amssymb,amsthm}
\usepackage{amsmath}
\usepackage{amscd}
\usepackage{xy}
\xyoption{all}
\usepackage[dvips]{graphicx}

\newtheorem{thm}{Theorem}

\begin{document}

\setlength{\baselineskip}{16pt}

\date{}

\title{An identity in Rota-Baxter algebras}

\author{ Rafael D\'\i az and  Marcelo P\'aez}

\maketitle

\begin{abstract}
We give explicit formulae and study the combinatorics of an
identity holding in all Rota-Baxter algebras. We describe the
specialization of this identity for a couple of examples of
Rota-Baxter algebras.
\end{abstract}


\section{Introduction}

The study of Rota-Baxter algebras was initiated by Baxter in his
works \cite{bax1, bax2}. The theory was later taken over by Rota
\cite{rota}, who gave an explicit construction of free Rota-Baxter
algebras and uncovered the relationship with symmetric functions.
Soon after Cartier studied free Rota-Baxter algebras in
\cite{cartier}. In the last few years the theory of Rota-Baxter algebras has received a
great impulse, mainly because of its applications to
renormalization, as formalized by Connes and Kreimer in
\cite{ck1, ck2, ck3}.  New techniques and
applications of Rota-Baxter algebras have been found by an active
group of researches in a number of important works, among which we
cite just a few \cite{a, guo, guo2, pay}.\\

In this work we consider the seemingly naive problem of writing an
element of the form $P^{a}(x)P^{b}(y)$ in a Rota-Baxter algebra as
a linear combination of terms of the form $P^{j}(xP^{i}(y))$ and
$P^{j}(P^{i}(x)y)$ with $i$ and $j$ varying. The existence of such
a linear combination is an immediate consequence of the
Rota-Baxter identity satisfied by the operator $P$. The actual
problem is to determine the coefficients involved in such an
expression as explicitly as possible. We provide a solution to
this problem, and take a look at its meaning in a couple of
Rota-Baxter algebras. We approach our problem from a rather
pedestrian point of view using a graphical notation to illustrate
our ideas.

\section{Basic ideas}
Let us fix $k$ a field of characteristic zero. A Rota-Baxter
algebra is a triple $(A,\lambda,P)$ where $A$ is an associative
$k$-algebra, $\lambda$ is a constant in $k$, and $P:A
\longrightarrow A$ is a $k$-linear operator satisfying the identity
$$P(x)P(y)=P(xP(y))+P(P(x)y)+\lambda P(xy)$$
for $x,y \in A$.  We find it convenient to use a graphical
notation to express our results.
\newpage
We represent the product on $A$ by
\begin{figure}[h!]
\begin{center}
\includegraphics[height=1.2cm]{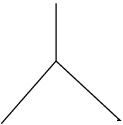}
\end{center}
\end{figure}

and the Rota-Baxter operator by

\begin{figure}[h!]
\begin{center}
\includegraphics[height=1.2cm]{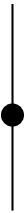}
\end{center}
\end{figure}

The Rota-Baxter identity satisfied by $P$ is represented
graphically by

\begin{figure}[h!]
\begin{center}
\includegraphics[height=1.2cm]{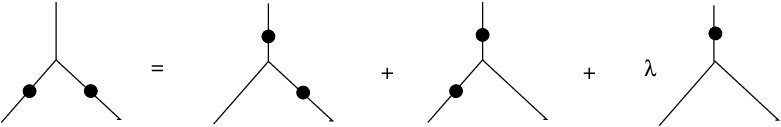}
\end{center}
\end{figure}

For example using the graphical form of the Rota-Baxter identity
one can see that

\begin{figure}[h!]
\begin{center}
\includegraphics[height=1.2cm]{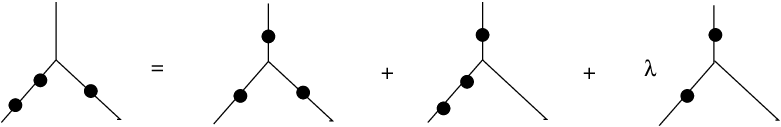}
\end{center}
\end{figure}

A further application of the graphical Rota-Baxter identity yields

\begin{figure}[h!]
\begin{center}
\includegraphics[height=1.2cm]{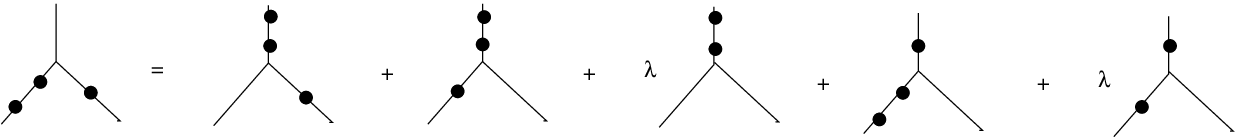}
\end{center}
\end{figure}
Thus we have shown that the following identity holds in any
Rota-Baxter algebra:
$$P^{2}(x)P(y)= P^{2}(xP(y)) + P^{2}(P(x)y) + \lambda P^{2}(xy) + P(P^{2}(x)y) + \lambda P(P(x)y).$$

The symbol $T(a,b,c)$ has two different meanings in this work:
\begin{itemize}
\item On
the one hand it stands for the operator $P^{c}(m \circ (P^{a}
\otimes P^{b})): A \otimes A \longrightarrow A$ where $m$ denotes
the product on $A$.

\item On the other hand it represents the tree with $a$ dots on the left
leg, $b$ dots on the right leg, and $c$ dots on the neck. The
three $T(a,b,c)$ is drawn as follows

\begin{figure}[h!]
\begin{center}
\includegraphics[height=1.2cm]{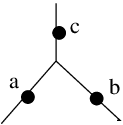}
\end{center}
\end{figure}
\end{itemize}

For example the three $T(1,2,3)$ is represented graphically as
follows:

\begin{figure}[h!]
\begin{center}
\includegraphics[height=1.2cm]{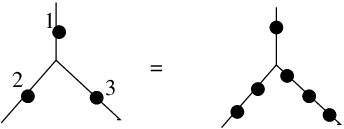}
\end{center}
\end{figure}

\newpage
It is clear from the graphical Rota-Baxter identity that each tree

\begin{figure}[h!]
\begin{center}
\includegraphics[height=1.2cm]{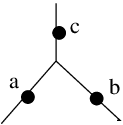}
\end{center}
\end{figure}

can be written as a linear combination with coefficients in
$k[\lambda]$ of trees of the form

\begin{figure}[h!]
\begin{center}
\includegraphics[height=1.2cm]{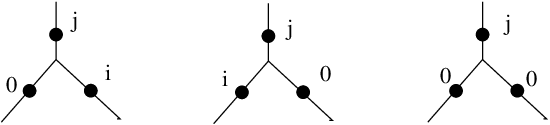}
\end{center}
\end{figure}
Indeed each application of the graphical Rota-Baxter relation
replaces the tree

\begin{figure}[h!]
\begin{center}
\includegraphics[height=1.2cm]{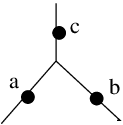}
\end{center}
\end{figure}

by the sum of trees

\begin{figure}[h!]
\begin{center}
\includegraphics[height=1.2cm]{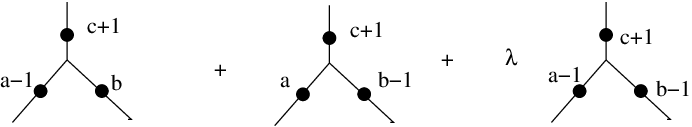}
\end{center}
\end{figure}

From an algorithmic point of view the graphical Rota-Baxter
identity can be described as the application of three weighted
moves:

\begin{enumerate}
\item{A weight $1$ move where a dot from the left leg moves up.  }
\item{A weight $1$ move where a dot from the right leg moves up.}
\item{A weight $\lambda$ move where a couple of dots, one from the left leg and another one from the
right leg, merge and move up as one dot.}

\end{enumerate}

We are ready to formulate our main results.

\section{Restricted case $\lambda =0$}

The case of Rota-Baxter algebras with $\lambda=0$ simplifies
considerably. We report on this special case because of its
applications and elegant proof.

\begin{thm}\label{HI}{\em
Let $a,b > 1$ and $c \geq 0$ be integers. The following identity
holds in any Rota-Baxter algebra
\[T(a,b,c)=\sum_{i=1}^{b}\binom{a-1+b-i}{ a-1
}T(0,i,a+b+c-i) + \sum_{i=1}^{a}\binom{b-1+a-i}{ b-1
}T(i,0,a+b+c-i).
\]}
\end{thm}

\begin{proof}
We justify only the left summand, the right summand is justified
in an analogous way. For $\lambda = 0$ only move 1 and move 2 are
allowed. With each move a dot from one of the legs moves up.
Suppose that after applying several times the Rota-Baxter identity
to $T(a,b,c)$ we arrive at a tree of the form $T(0,i,j)$. Then a
total of $a + b - i$ dots from the legs have moved up so $j=
a+b+c-i$. Necessarily the last dot moving up comes from the left
leg. The other dots moved up in an arbitrary order, so this
explain the factor $$\binom{a-1+b-i}{ a-1 }.$$
\end{proof}

Consider the Rota-Baxter algebra $(C(\mathbb{R}),0 ,P),$ where
$C(\mathbb{R})$ denotes the algebra of continuous functions on
$\mathbb{R}$ and $P$ is the Riemann integral operator given by
$$P(f)(y)=\int_{0}^{y}f(x)dx.$$ For real numbers $0 \leq x \leq y$ and $a \in  \mathbb{N}_{+},$
let $\Delta_{x,a}^{y}$ be the convex polytope
$$\Delta_{x,a}^{y} = \{(x_1,...,x_a) \in (\mathbb{R}_{\geq 0})^n \ \ | \ \ x\leq x_1 \leq...\leq x_a \leq y \}.$$
For $a \geq 1$ we let $v_a(x,y)$ be the volume of
$\Delta_{x,a}^{y}.$ By convention we set $v_0(x,y)=1.$  It is easy
to check that
$$P^{a+1}(f)(y)=\int_{0}^{y}f(x)v_a(x,y)dx.$$ Theorem \ref{HI} implies the following
result.

\begin{thm}
$$\int_{0}^{z}\int_{0}^{z}f(x)g(y)v_a(x,z)v_b(y,z)dxdy=$$
\begin{align*}
&&\sum_{i=1}^{b+1}\binom{a+b+1-i}{a}\int_{0 \leq x \leq y \leq
z}g(x)f(y)v_{i-1}(x,y)v_{a+b+1-i}(y,z)dxdy\\
&+&\sum_{i=1}^{a+1}\binom{b+a+1-i}{b}
\int_{0 \leq x \leq y \leq z}f(x)g(y)v_{i-1}(x,y)v_{b+a+1-i}(y,z)dxdy.
\end{align*}
\end{thm}
\section{Generic case}

Now we considerer the generic situation, i.e., a Rota-Baxter
algebra with $\lambda \neq 0$.
\begin{thm}\label{FI}{\em
Let $a,b \geq 1$ and $c \geq 0$ be integers. The following
identity holds in any Rota-Baxter algebra:
\begin{align*}
T(a,b,c)&=\sum_{(i,j) \in D_1}c_1(a,b;i,j)T(0,i,c+j) \\
&+\sum_{(i,j) \in D_2}c_2(a,b;i,j)T(0,i,c+j)\\
&+\sum_{(i,j) \in D_3}c_3(a,b;i,j)T(i,0,c+j)\\
&+\sum_{(i,j) \in D_4}c_4(a,b;i,j)T(i,0,c+j)\\
&+\sum_{j \in D_5}c_5(a,b;j)T(0,0,c+j)
\end{align*}
where
\begin{align*}
D_1&=\{(i,j)\in \mathbb{N}_+ \times \mathbb{N}_+ \ \ |\ \  1 \leq i \leq b,\ \  b-i+1 \leq j, \ \ a \leq j,\ \ j \leq a +b -i \}, \\
D_2&=\{(i,j)\in \mathbb{N}_+ \times \mathbb{N}_+ \ \ |\ \  1 \leq i \leq b-1,\ \  b-i \leq j, \ \  a \leq j,\ \ j \leq a +b-i -1 \},\\
D_3&=\{(i,j)\in \mathbb{N}_+ \times \mathbb{N}_+ \ \ |\ \  1 \leq i \leq a,\ \ b \leq j, \ \ a-i +1 \leq j,\ \ j \leq a +b -i \},\\
D_4&=\{(i,j)\in \mathbb{N}_+ \times \mathbb{N}_+ \ \ |\ \  1 \leq i \leq a-1,\ \ b \leq j, \ \ a-i \leq j,\ \ j \leq a +b-i -1 \},\\
D_5&=\{ j \in \mathbb{N} \ \ |\ \ b
\leq j,\ \ a \leq j,\ \ j \leq a +b -1 \},
\end{align*}

and
\begin{align*}
c_1(a,b;i,j) &=\binom{j-1}{i+j-b-1,\ \ j-a,\ \ a+b-i-j}\lambda^{a+b-i-j}, \\
c_2(a,b;i,j)&=\binom{j-1}{i+j-b,\ \ j-a, \ \ a+b-i-j-1}\lambda^{a+b-i-j},\\
c_3(a,b;i,j)&=\binom{j-1}{j-b, \ \ i+j-a-1,\ \ a+b-i-j}\lambda^{a+b-i-j},\\
c_4(a,b;i,j)&=\binom{j-1}{j-b, \ \ i+j-a, \ \ a+b-i-j-1}\lambda^{a+b-i-j},\\
c_5(a,b;j)&=\binom{j-1}{ j-b,\ \ j-a, \ \ a+b-j-1}\lambda^{a+b-j}.
\end{align*}}
\end{thm}

\begin{proof}
Notice that with each move a dot is added to the neck. If starting
from the tree $T(a, b, c)$ we arrive using the allowed moves to
the graph $T(0,i,c+j)$, then necessarily we must have applied $j$
moves and the last move must have been either move $1$ or move
$3$. Let us consider the case were the last move is of type $1$.
The other $j-1$ moves are distributed into $k_1$ moves of type
$1$, $k_2$ moves of type $2$, and $k_3$ moves of type $3$,  given
rise to the combinatorial number
$$\binom{j-1}{ k_1, k_2, k_3}.$$
The numbers $k_1,k_2$ and $k_3$ are subject to the constrains
$$k_1 + k_2 + k_3 = j-1, \ \  k_1 + k_3 = a - 1 \mbox { and } k_2 + k_3 = b - i.$$ Solving this linear system of equations we find that
$$k_1 = i+j-b-1, \ \ k_2 =j-a \mbox { and } k_3 = a+b-i-j.$$
This justifies the expression for $c_1(a,b,c;i,j)$ from the
statement of the Theorem. We proceed to justify the expression for
$c_2(a,b,c;i,j)$ which arises when the last move taken in the path
towards $T(0,i,j)$ is of type $3$. The remaining new $j-1$ dots in
the neck move up as consequence of the application of any of the
moves, giving rise to the factor $\binom{j-1}{k_1, k_2, k_3}$
where $k_1, k_2$ and $k_3$ are subject to the constraints
$$k_1 + k_2 + k_3 = j-1, \ \ k_1 + k_3 = a-1, \ \ k_2 + k_3 = b - i -1 .$$
Solving this equations we find that
$$k_1 = i+j-b, \ \ k_2 = j-a, \ \ k_3 =a+b-i-j-1,$$
thus we have justified the factor $$\binom{j-1}{i+j-b, j-a,
a+b-i-j-1}\lambda^{a+b-i-j}$$ appearing in the formula for
$c_2(a,b,c;i,j)$. The formulae for $c_3(a,b,c;i,j)$ and
$c_4(a,b,c;i,j)$ are derived in a fairly similar way. Let us
consider the formula for $c_5(a,b,c;j).$ In this case the last
move is necessarily of type $3$ and gives rise to the factor
$\binom{j-1}{k_1, k_2, k_3}$ where $k_1, k_2$ and $k_3$ satisfy
the constrains
$$k_1 + k_2 + k_3 = j-1, \ \ k_1 + k_3 = a-1, \ \ k_2 + k_3 = b-1 .$$
We find that
$$k_1 = j-b, \ \ k_2 = j-a, \ \ k_3 =a+b-j-1,$$
which justifies the factor $$\binom{j-1}{j-b, j-a,
a+b-j-1}\lambda^{a+b-j}$$ appearing in the formula for
$c_5(a,b;j)$.

\end{proof}

Consider the Rota-Baxter algebra $(\{f:
\mathbb{N}
\longrightarrow \mathbb{C}\},-1, P)$, where the operator $$P:\{f:
\mathbb{N}
\longrightarrow \mathbb{C}\} \longrightarrow \{f:
\mathbb{N}
\longrightarrow \mathbb{C}\}$$ is given by the Riemann sum $$P(f)(m)= \sum_{n=1}^{
m}f(n).$$ For $a \geq 1$ one can check that
$$P^{a}f(m) = \sum_{1\leq n_1 \leq ... \leq n_{a} \leq m}f(n_1).$$
In particular
$$P^{a}1(m) = |\Omega_{a}^{m}|,$$ where $$\Omega_{a}^{m}= \{(n_1,...,n_{a}) \in \mathbb{N}_+ \ \ | \ \
1\leq n_1 \leq ... \leq n_{a} \leq m \}.$$\\
It is not hard to show using the Chu-Vandermonde identity that
$$|\Omega_{a}^{m}|= \sum_{s=1}^{m}\binom{m}{s}\binom{a}{s}=\binom{a+m}{m} -1.$$ As a consequence of Theorem \ref{FI} we get that the numbers
$|\Omega_{a}^{m}|$ satisfy the following identity.

\begin{thm}{\em For integers $a,b \geq 1$ we have
\begin{align*}
|\Omega_{a}^{m}||\Omega_{b}^{m}|&=\sum_{(i,j) \in D_1}[c_1(a,b;i,j)(-1)]|\Omega_{i+j}^{m}| \\
&+\sum_{(i,j) \in D_2}[c_2(a,b;i,j)(-1)]|\Omega_{i+j}^{m}|\\
&+\sum_{(i,j) \in D_3}[c_2(a,b;i,j)(-1)]|\Omega_{i+j}^{m}|\\
&+\sum_{(i,j) \in D_4}[c_4(a,b;i,j)(-1)]|\Omega_{i+j}^{m}|\\
&+\sum_{j \in D_5}[c_5(a,b;j)(-1)]|\Omega_{j}^{m}|.
\end{align*}}
\end{thm}

\section{Acknowledgment}

We thank Dominique Manchon for helpful comments and suggestions.
Thanks also to Takasi Kimura, Jos\'e Mijares, Eddy Pariguan and
Sylvie Paycha.

\vspace{8mm}

\noindent ragadiaz@gmail.com, \ \ mathmarcelo@gmail.com.\\
\noindent Escuela de Matem\'aticas, Universidad Central de Venezuela,
Caracas 1020, Venezuela.\\

\end{document}